\documentclass[11pt]{amsart}

\usepackage{latexsym,amsfonts,amssymb}
\usepackage{amsmath,amsthm,amscd}
\usepackage{hyperref,psfrag}

\usepackage[dvips]{epsfig}
\usepackage{psfrag}

\usepackage{graphicx}

\usepackage{a4wide}

\normalfont\upshape

\newcommand{\sign}{\text{sign}}

\usepackage{fancyheadings}
\usepackage[all]{xy}
\SelectTips{cm}{}

\theoremstyle{definition}
\newtheorem{thm}{Theorem}[section]
\newtheorem*{unthm}{Theorem}
\newtheorem{cor}[thm]{Corollary}

\newtheorem{lem}[thm]{Lemma}
\newtheorem{rem}[thm]{Remark}
\newtheorem{prop}[thm]{Proposition}
\newtheorem{defn}[thm]{Definition}
\newtheorem{example}[thm]{Example}

\numberwithin{equation}{section}

\usepackage{bbm}
\def\C{{\mathbbm C}}

\def\R{{\mathbbm R}}
\def\Z{{\mathbbm Z}}

\def\ZPl{\mathbb{Z}Pl}
\def\xt{\widetilde{x}}

\renewcommand{\to}{\rightarrow}

\def\sh{\mathrm{sh}}
\def\ct{\mathrm{ct}}
\def\opol{\mathrm{OPol}}
\def\osym{\mathrm{O}\Lambda}
\def\SSYT{SSYT}
\def\Gr{\mathrm{Gr}}
\def\P{\mathbb{P}}

\def\Ind{{\mathrm{Ind}}}

\def\shuffle{\,\raise 1pt\hbox{$\scriptscriptstyle\cup{\mskip
               -4mu}\cup$}\,}

\newcommand{\refequal}[1]{\xy {\ar@{=}^{#1}
(-1,0)*{};(1,0)*{}};
\endxy}

\usepackage[vcentermath,enableskew]{youngtab}

\title{The odd Littlewood-Richardson rule}
      \author{Alexander P. Ellis}

\date{November 15, 2011}

\newcommand{\sym}{\mathrm{\Lambda}}

\begin{document}
%

\maketitle

\begin{abstract}
In previous work with Mikhail Khovanov and Aaron Lauda we introduced two odd analogues of the Schur functions: one via the combinatorics of Young tableaux (odd Kostka numbers) and one via the odd symmetrization operator.  In this paper we introduce a third analogue, the plactic Schur functions.  We show they coincide with both previously defined types of Schur function, confirming a conjecture.  Using the plactic definition, we establish an odd Littlewood-Richardson rule.  We also re-cast this rule in the language of polytopes, via the Knutson-Tao hive model.
\end{abstract}

\setcounter{tocdepth}{1} \tableofcontents


%
\section{Introduction}
%

%
\subsection{The classical story}\label{subsec-classical-story}
%

Schur functions play a central role in the beautiful circle of ideas around symmetric functions, symmetric and general linear group representations, and linear enumerative geometry.  On the one hand, the Schur functions $s_\lambda$ corresponding to partitions $\lambda$ of $k$ give an integral orthonormal basis of the algebra of symmetric functions $\sym$ in degree $k$.  Under the Frobenius correspondence (see, for instance, \cite[Section 7.18]{Stanley2})
\begin{equation*}
\bigoplus_{n\geq0}K_0(\C[S_n])\cong\sym,
\end{equation*}
the Schur function $s_\lambda$ corresponds to the irreducible representation $L_\lambda$ associated to the partition $\lambda$. Recall that the product on the symmetric group side is given as follows.  If $V$ is a representation of $S_k$ and $W$ is a representation of $S_\ell$, then
\begin{equation}
[V]\cdot[W]=[\Ind_{S_k\times S_\ell}^{S_{k+\ell}}(V\otimes W)].
\end{equation}
On the other hand, the ring $\sym_n$ obtained by setting $s_{(1^{n+1})}=s_{(1^{n+2})}=\ldots=0$ can be viewed as the character ring for polynomial representations of $GL_n(\C)$.  Under this correspondence, the image of $s_\lambda$ in $\sym_n$ equals the character of the irreducible representation $V_\lambda$ associated to $\lambda$.  This time the product is the usual product of characters, that is, $s_\mu s_\nu$ is the character of $V_\mu\otimes V_\nu$.  Either of these two descriptions immediately implies that the structure coefficients of a product
\begin{equation}\label{eqn-intro-lrc}
s_\mu s_\nu=\sum_\lambda c_{\mu\nu}^\lambda s_\lambda
\end{equation}
are all non-negative integers; they are called \textit{Littlewood-Richardson coefficients}.  The determination of these coefficients is a classical and important problem.  By the above, $c_{\mu\nu}^\lambda$ has interpretations as
\begin{enumerate}
\item the multiplicity of $L_\lambda$ in $\Ind_{S_k\times S_\ell}^{S_{k\times\ell}}(L_\mu\otimes L_\nu)$,
\item the multiplicity of $V_\lambda$ in $V_\mu\otimes V_\nu$, and
\item the dimension of the space of $GL_n(\C)$-invariant vectors in $V_\mu\otimes V_\nu\otimes V_\lambda^*$.
\end{enumerate}
\noindent The combinatorics of Young tableaux become an important tool at this point: partitions correspond to Young diagrams, and the coefficients $c_{\mu\nu}^\lambda$ can be computed by counting certain skew tableaux.  We review this Littlewood-Richardson rule in Section \ref{subsec-even-lr}; more complete expositions can be found in \cite[Chapter 5]{Fulton} and \cite[Appendix A1]{Stanley2}.

Schur functions also arise in linear enumerative geometry.  The cohomology ring of the Grassmannian $\Gr(k,n)$ of complex $k$-planes in $\C^n$ is isomorphic to the quotient of $\sym_k$ by the ideal generated by all complete symmetric functions $h_m$ for $m>n-k$.  Geometrically, elementary functions $e_m$ correspond to the Chern classes of the tautological $k$-plane bundle and complete functions $h_m$ correspond (up to sign) to the Chern classes of a complementary $(n-k)$-plane bundle.  Under this identification, the image of the Schur function $s_\lambda$ in $H^*(\Gr(k,n))$ equals zero unless $\lambda_1\leq n-k$ and $\ell(\lambda)\leq k$.  The nonzero Schur functions give a basis of the cohomology ring, and $s_\lambda$ is Poincar\'{e} dual to the class of a corresponding \textit{Schubert variety}.  Relative to a fixed full flag $0=V_0\subset V_1\cdots\subset V_n=\C^n$, a point $W\in\Gr(k,n)$ is in the Schubert variety $\overline{C_\lambda}$ for the partition $\lambda$ if and only if
\begin{equation*}
\dim(W\cap V_i)\geq i+\lambda_i-(n-k)
\end{equation*}
for $1\leq i\leq n$.  By \eqref{eqn-intro-lrc}, then, $c_{\mu\nu}^\lambda$ is the coefficient of $[\overline{C_\lambda}]$ in the cup product $[\overline{C_\mu}]\cdot[\overline{C_\nu}]$.  Therefore we have a geometric interpretation of $c_{\mu\nu}^\lambda$ as
\begin{enumerate}
\item[(4)] the multiplicity of $[\overline{C_\lambda}]$ in the cup product $[\overline{C_\mu}]\cdot[\overline{C_\nu}]$.
\end{enumerate}
When $|\mu|+|\nu|=\dim(\Gr(k,n))=k(n-k)$, the coefficient $c_{\mu\nu}^{(n-k)^k}$ is the (finite!) number of $k$-planes satisfying the dimension criteria of both $\overline{C_\mu}$ and $\overline{C_\nu}$.  For example, one can use this analysis to compute the number of lines which intersect five general 3-planes in $\P^6=\P(\C^7)$.

All of the above is to say that Schur functions provide important insight into areas of geometry, algebra, and representation theory.  The goal of this paper is to provide a small first step in this direction in the \textit{odd} setting, in a sense we will now describe.

%
\subsection{Outline of this paper}
%

An ongoing project initiated in joint work with Khovanov and Lauda is an attempt to give a construction of the Ozsv\'{a}th-Rasmussen-Szab\'{o} odd Khovanov homology \cite{ORS} via a 2-representation-theoretic approach analogous to Webster's construction of (even) Khovanov homology \cite{Web}, \cite{Web2}.  As a byproduct of this investigation, an ``odd'' analogue $\osym$ of the algebra $\sym$ was defined and found to admit signed analogues of many of the combinatorial properties of $\sym$ \cite{EK}, \cite{EKL}.  We gave two different candidates for elements playing the role analogous to that of Schur functions.  The present work gives a third candidate, proves all three definitions are equivalent, and proves an odd analogue of the Littlewood-Richardson rule.  

Section \ref{subsec-odd-symmetric} reviews the definitions and basic properties of odd symmetric functions and odd Schur functions; the short Section \ref{subsec-boxterpretations} gives a convenient notation for the many signs which arise in the odd setting.  Section \ref{subsec-odd-plactic} contains the new definition of odd Schur functions, and Section \ref{subsec-comparison} proves the following.
\begin{unthm} The three definitions of odd Schur functions all coincide.\end{unthm}
\noindent In Section \ref{subsec-even-lr} we review the even Littlewood-Richardson rule, and in Section \ref{subsec-odd-lr} we formulate and prove an odd analogue.
\begin{unthm} The odd Littlewood-Richardson coefficient $c_{\mu\nu}^\lambda$ is a signed count of semistandard skew tableaux $S$ of shape $\lambda/\mu$ and content $\nu$ such that the row word of $S$ is Yamanouchi.\end{unthm}
\noindent Finally, in Section \ref{subsec-hives}, we re-cast the odd Littlewood-Richardson rule in the language of Knutson-Tao hives.

%
\subsection{Acknowledgments}
%

The author thanks Mikhail Khovanov for many helpful conversations and suggestions.  For the duration of this work, the author was supported by an NSF Graduate Research Fellowship.

%
\section{Review of odd symmetric functions}
%

%
\subsection{Odd symmetric polynomials}\label{subsec-odd-symmetric}
%

Now we will briefly review the constructions of \cite{EK}, \cite{EKL}.  Let $n\geq1$ and let
\begin{equation*}
\opol_n=\Z\langle x_1,\ldots,x_n\rangle/(x_ix_j+x_jx_i\text{ if }i\neq j)
\end{equation*}
be the ring of skew polynomials in $n$ variables.  The \textit{ascend-sorting} of a monomial $x_{i_1}\cdots x_{i_r}$ in $\opol_n$ is defined to be the monomial obtained by sorting the subscripts $i_1,\ldots,i_r$ into non-decreasing order without introducing any sign.  For example, the ascend-sorting of $x_1x_3x_2$ is $x_1x_2x_3$.  We will sometimes write $:X:$ for the ascend-sorting of a monomial $X$.

Although the generators $x_1,\ldots,x_n$ pairwise supercommute, the algebra $\opol_n$ is not supercommutative:
\begin{equation*}
(x_1+x_2)x_1=x_1(x_1-x_2).
\end{equation*}
The algebra $\opol_n$ is finite dimensional in each degree and has no zero divisors, but it does not have unique factorization:
\begin{equation*}
(x_1-x_2)^2=x_1^2+x_2^2=(x_1+x_2)^2.
\end{equation*}

Have the symmetric group $S_n$ act on the free algebra $\Z\langle x_1,\ldots,x_n\rangle$  as 
\begin{equation}\label{eqn-symm-action}
s_i(x_j)=\begin{cases}-x_{i+1}&\text{if }j=i,\\-x_i&\text{if }j=i+1,\\-x_j&\text{if }j\neq i,i+1,\end{cases}
\end{equation}
where $s_i$ is the transposition of $i$ and $i+1$ and the action of $s_i$ is a ring endomorphism.  For $i=1,\ldots,n-1$, the $i$-th \textit{odd divided difference operator} is the $\Z$-linear map $\partial_i:\Z\langle x_1,\ldots,x_n\rangle\to\Z\langle x_1,\ldots,x_n\rangle$ defined by
\begin{equation}\begin{split}
&\partial_i(x_j)=\begin{cases}1&\text{if }j=i,i+1,\\0&\text{otherwise,}\end{cases}\\
&\partial_i(fg)=\partial_i(f)g+s_i(f)\partial_i(g).
\end{split}\end{equation}
The second line is called the Leibniz rule.  It is easy to check that
\begin{equation*}
\partial_i(x_jx_k+x_kx_j)=0
\end{equation*}
for all $j,k$, so we can consider $\partial_i$ as an operator on $\opol_n$.  On $\opol_n$, the kernel and image of $\partial_i$ coincide; by contrast with the even case, however, these do not equal the space of invariants or anti-invariants of the action of $s_i$.  

\begin{defn}The ring of \textit{odd symmetric polynomials} is the subring
\begin{equation}
\osym_n=\bigcap_{i=1}^{n-1}\ker(\partial_i)\subset\opol_n.
\end{equation}\end{defn}

We endow $\osym_n$ with a $\Z$- and a $\Z/2$-grading, whose degree functions we denote by $\deg_\Z$ and $\deg_s$ respectively.  Write $\deg(f)=(\deg_\Z(f),\deg_s(f))$ for short.  The degree of each $x_i$ is defined to be
\begin{equation*}
\deg(x_i)=(2,1).
\end{equation*}
We will consider $\osym_n$ as a superalgebra via the $\Z/2$-grading.  In particular, the product structure on $\osym_n\otimes\osym_n$ is, on homogeneous elements,
\begin{equation*}
(f\otimes g)(f'\otimes g')=(-1)^{\deg_s(g)\deg_s(f')}(ff')\otimes(gg').
\end{equation*}

\vspace{0.07in}

For $i=1,\ldots,n$, let $\xt_i=(-1)^{i-1}x_i$.  We define elements of $\opol_n$ for each $k\geq1$,
\begin{equation*}\begin{split}
&e_k=\sum_{1\leq i_1<\ldots<i_k\leq n}\xt_{i_1}\cdots\xt_{i_k},\\
&h_k=\sum_{1\leq i_1\leq\ldots\leq i_k\leq n}\xt_{i_1}\cdots\xt_{i_k}.
\end{split}\end{equation*}
Set $e_0=h_0=1$ and $e_k=h_k=0$ for $k<0$.  If $k>n$, then $e_k=0$.  Both these families of skew polynomials are odd symmetric, and they satisfy the relation
\begin{equation*}
\sum_{k=0}^\ell(-1)^{\frac{1}{2}k(k+1)}e_kh_{\ell-k}=0\qquad\text{if }\ell\geq1.
\end{equation*}
They also satisfy the \textit{odd defining relations},
\begin{equation}\label{eqn-ODRs}\begin{split}
&e_ae_b=e_be_a\quad\text{if }a+b\text{ is even},\\
&e_ae_b+(-1)^ae_be_a=(-1)^ae_{a+1}e_{b-1}+e_{b-1}e_{a+1}\quad\text{if }a+b\text{ if odd}.
\end{split}\end{equation}
The relations \eqref{eqn-ODRs} also hold if all $e$'s are replaced by $h$'s.

This next proposition follows from Section 2.1 of \cite{EKL}.
\begin{prop}\label{prop-odd-properties} The algebra $\osym$ has a presentation by generators $e_0=1,e_1,e_2,\ldots$ and relations \eqref{eqn-ODRs}.  A basis for $\osym$ in $\Z$-degree $2k$ is given by all products $e_{\lambda_1}\cdots e_{\lambda_r}$  with $\lambda_1\geq\ldots\lambda_r\geq1$, $|\lambda_1|+\ldots+|\lambda_r|=k$.  The same is true if all $e$'s are replaced by $h$'s.\end{prop}
The products $e_\lambda=e_{\lambda_1}\cdots e_{\lambda_r}$ and $h_\lambda=h_{\lambda_1}\cdots h_{\lambda_r}$ with $\lambda_1\geq\ldots\geq\lambda_r$ are called \textit{odd elementary symmetric functions} and \textit{odd complete symmetric functions}, respectively.  The sets of elementary and complete functions in $\Z$-degree $2k$ are naturally indexed by all partitions of $k$.

\vspace{0.07in}

There are maps
\begin{equation*}\begin{split}
&\opol_{n+1}\to\opol_n\\
&x_i\mapsto x_i\text{ if }1\leq i\leq n,\\
&x_{n+1}\mapsto0.
\end{split}\end{equation*}
These induce maps $\osym_{n+1}\to\osym_n$ which send $h_k\mapsto h_k$ and $e_k\mapsto e_k$ for all $k$.  The inverse limit in the category of graded rings of the resulting system is called the ring of \textit{odd symmetric functions},
\begin{equation*}
\osym=\underleftarrow{\lim}\osym_n.
\end{equation*}
The ring $\osym$ can be given the structure of a Hopf superalgebra \cite{EK}.  The coproduct is
\begin{equation*}
\Delta(h_k)=\sum_{i+j=k}h_i\otimes h_j,\qquad\Delta(e_k)=\sum_{i+j=k}e_i\otimes e_j.
\end{equation*}
As a Hopf superalgebra, $\osym$ is neither commutative nor cocommutative.

\vspace{0.07in}

As a result of Proposition \ref{prop-odd-properties}, we can define certain symmetries of $\osym$,
\begin{equation*}\begin{split}
&\psi_1(h_k)=e_k\\
&\qquad\text{bialgebra automorphism (not an involution),}\\
&\psi_2(h_k)=(-1)^{\frac{1}{2}k(k+1)}h_k\\
&\qquad\text{algebra involution (not a coalgebra homomorphism),}\\
&\psi_3(h_k)=h_k\\
&\qquad\text{algebra anti-involution (not a coalgebra homomorphism).}
\end{split}\end{equation*}
The map $\psi_1\psi_2$ is an involution, and both $\psi_1,\psi_2$ commute with $\psi_3$.  The antipode of the Hopf superalgebra structure on $\osym$ is $S=\psi_1\psi_2\psi_3$.  For a partition $\lambda=(\lambda_1,\ldots,\lambda_r)$, define
\begin{equation*}
\epsilon_\lambda=(-1)^{\sum_j\frac{1}{2}\lambda^T_j(\lambda^T_j-1)},
\end{equation*}
where $\lambda^T$ is the transpose (or dual, or conjugate) partition to $\lambda$.  That is, $\lambda^T_j$ is the height of the $j$-th column of $\lambda$.  The involution $\psi_1\psi_2$ exchanges the elementary and complete bases up to sign,
\begin{equation}\label{eqn-psi12-eh}
\psi_1\psi_2(e_\lambda)=(-1)^{|\lambda|}\epsilon_{\lambda^T}h_\lambda,\qquad\psi_1\psi_2(h_\lambda)=(-1)^{|\lambda|}\epsilon_{\lambda^T}e_\lambda.
\end{equation}

\vspace{0.07in}

In \cite{EK}, an odd analogue of the Schur functions was given as follows.  From here on, we will identify partitions and Young diagrams without comment.  Recall that a \textit{Young tableau} $T$ on a partition $\lambda$ with entries in an ordered alphabet $A$ is an assignment of an element of $A$ to each box of $\lambda$.  When considering $\osym$ we will take $A=\Z_{>0}$ and when considering $\osym_n$ we will take $A=\lbrace1,\ldots,n\rbrace$.  A tableau $T$ is called \textit{semistandard} if its entries are non-decreasing in rows (left to right) and strictly increasing in columns (top to bottom).  The \textit{row word} $w_r(T)$ of a tableau $T$ is the word in the alphabet $A$ obtained by reading the entries of $T$ from left to right, bottom to top.  For example,
\begin{equation*}
T=\young(112,23)\text{ is a semistandard Young tableau of shape }(3,2),\text{ and }w_r(T)=23112.
\end{equation*}
The \textit{content} of a tableau $T$ is the tuple $(a_1,\ldots,a_r)$, where $a_i$ is the number of entries of $T$ equal to $i$.  For each partition $\lambda$, there is a unique Young tableau $T_\lambda$ of shape and content both equal to $\lambda$.  For example,
\begin{equation*}
T_{(21)}=\young(11,2),\qquad T_{(311)}=\young(111,2,3).
\end{equation*}
If $T$ is a Young tableau, write $\sh(T)$ for its shape and $\ct(T)$ for its content.  Let $\SSYT(\lambda)$ be the set of semistandard Young tableaux of shape $\lambda$ and let $\SSYT(\lambda,\mu)$ be the set of semistandard Young tableaux of shape $\lambda$ and content $\mu$.

In order to use tableaux in the odd setting, define the sign of a tableau $T$, denoted $\sign(T)$, to be the sign of the minimal length permutation which sorts $w_r(T)$ into non-decreasing order.  For instance, with $T$ as above,
\begin{equation*}
\sign(T)=(-1)^5=-1,\qquad\sign(T_{(21)})=(-1)^2=1,\qquad\sign(T_{(311)})=(-1)^7=-1.
\end{equation*}
For partitions $\lambda,\mu$, the \textit{odd Kostka number} $K_{\lambda\mu}$ is defined to be a signed count of semistandard tableaux of shape $\lambda$ and content $\mu$,
\begin{equation}
K_{\lambda\mu}=\sign(T_\lambda)\sum_{T\in\SSYT(\lambda,\mu)}\sign(T).
\end{equation}
It is readily verified that as a matrix, $(K_{\lambda\mu})_{|\lambda|=|\mu|=k}$ is lower triangular and unimodular when the partitions basis is ordered lexicographically, so we can define the family of \textit{combinatorial (or Kostka) odd Schur functions} by the change of basis relation
\begin{equation}
h_\mu=\sum_{\lambda\vdash k}K_{\lambda\mu}s^K_\lambda,
\end{equation}
where $\mu\vdash k$.

The following proposition combines several of the statements from Section 3.3 of \cite{EK}.
\begin{prop} The family $\lbrace s_\lambda^K\rbrace_{\lambda\vdash k}$ is an integral unimodular basis for $\osym$ in $\Z$-degree $2k$.  This family is signed-orthonormal,
\begin{equation}
(s_\lambda^K,s_\mu^K)=\epsilon_\lambda\delta_{\lambda\mu},
\end{equation}
with respect to the bilinear form of \cite{EK}.  The involution $\psi_1\psi_2$ and the anti-involution $\psi_3$ act on $s_\lambda$ as
\begin{equation}\label{eqn-psi-schur}
\psi_1\psi_2(s_\lambda^K)=(-1)^{\ell(w_\lambda)+|\lambda|}s_{\lambda^T}^K,\qquad\psi_3(s_\lambda^K)=\epsilon_\lambda\sign(T_\lambda)s_\lambda^K.
\end{equation}
Here, $w_\lambda$ is the element of $S_k$ which is combinatorially defined in Proposition 2.14 of \cite{EK} (it is the minimal representative of the unique double coset in $S_{\lambda^T}\backslash S_k/S_\lambda$ which gives a nonzero homomorphism from $\Ind_{S_\lambda}^{S_k}(V_\text{triv})$ to $\Ind_{S_{\lambda^T}}^{S_k}(V_\text{sign})$) and $\ell$ is the Coxeter length function.\end{prop}

In \cite{EKL}, another family of Schur functions was introduced.  Let
\begin{equation*}
\partial_{w_0}=\partial_1(\partial_2\partial_1)\cdots(\partial_{n-1}\cdots\partial_1),
\end{equation*}
a particular choice of longest odd divided difference operator (odd divided difference operators corresponding to other choices of reduced expression could differ by a factor of $-1$).  For a skew polynomial $f$, let $f^{w_0}$ be the result of acting by $w_0$ on $f$ via the action \eqref{eqn-symm-action}.  The \textit{odd-symmetrized Schur functions} are defined by
\begin{equation}
s_\lambda^s=(-1)^{\binom{n}{3}}\left(\partial_{w_0}(x_1^{\lambda_1}\cdots x_n^{\lambda_n}x_1^{n-1}x_2^{n-2}\cdots x_{n-1})\right)^{w_0}.
\end{equation}
For the motivation behind this name, see Section 2.2 of \cite{EKL}.  The goal of Section \ref{subsec-comparison} is to prove that for all $\lambda$, $s_\lambda^s=s_\lambda^K$.

%
\subsection{Boxterpretations}\label{subsec-boxterpretations}
%

It will be convenient to have a systematic way to handle the signs which arise in the odd setting.  Most combinatorial quantities arising in these signs can be described by counting certain boxes in some Young diagram or tableau; we call these descriptions ``boxterpretations.''  Here are some which have already arisen:
\begin{itemize}
\item For a Young diagram $\lambda$, the sum $\sum_j\frac{1}{2}\lambda^T_j(\lambda^T_j-1)$ can be described as: for each box $B$, add the number of boxes directly above $B$.  This is the exponent of $-1$ used in defining the sign $\epsilon_\lambda$.
\item For a Young diagram $\lambda$, the quantity $\ell(w_\lambda)$ can be described as: for each box $B$, add the number of boxes above and to the right of $B$.
\item For a Young tableau $T$: for each box $B$, add the number of boxes above and with smaller entry.  Then $\sign(T)$ equals $-1$ raised to this count.
\item As a special case of the previous, to obtain $\sign(T_\lambda)$, for each box $B$, count the number of boxes above $B$.
\end{itemize}
We will use notations adapted from \cite{Fulton} to help condense these descriptions.  If a box $B$ of a Young diagram is in a row above that of a box $B'$, we say $B$ is North of $B'$; if $B$ is North of or in the same row as $B'$, we say $B$ is north of $B'$.  Likewise for East/east, West/west, and South/south.  Write $N(B)$ for the number of boxes North of $B$, $sW(B)$ for the number of boxes southWest (both south and West) of $B$, and so forth.  If we are considering a Young tableau and decorate a direction with one of $\lbrace<,\leq,>,\geq\rbrace$, this means to only count boxes with entry lower than (lower than or equal to, greater than, greater than or equal to) $B$.  So $E^>(B)$ means the number of boxes East of $B$ with strictly greater entry than $B$.  If we evaluate one of these counting functions on a diagram or a tableau, we mean to sum over all boxes, evaluating the function at each.  For instance,
\begin{equation*}
SW^\geq(T)=\sum_{B\in T}\#\lbrace\text{boxes of }T\text{ SouthWest of }B\text{ with entry }\geq\text{ that of }B\rbrace.
\end{equation*}
One last decoration: $dN$ means directly North, that is, North and neither East nor West (and likewise for the other counts).

\begin{example} Let
\begin{equation*}
T=\young(1122,2334,344,56).
\end{equation*}
Then $dN(T)=16$, $E^>(T)=28$, and $sW(T)=47$.\end{example}

In this notation, we have the following boxterpretations:
\begin{equation*}
\epsilon_\lambda=(-1)^{dN(\lambda)},\qquad(-1)^{\ell(w_\lambda)}=(-1)^{NE(\lambda)},\qquad\sign(T)=(-1)^{N^<(T)},\qquad\sign(T_\lambda)=(-1)^{N(\lambda)}.
\end{equation*}

%
\section{Odd plactic Schur functions}
%

%
\subsection{Definition and basic properties}\label{subsec-odd-plactic}
%

Let $A=\lbrace a_1,a_2,\ldots\rbrace$ be an ordered alphabet.  In practice, we will take $A=\Z_{>0}$ when working with $\osym$ and $A=\lbrace1,2,\ldots,n\rbrace$ when working with $\osym_n$.  In order to add, multiply, and assign signs to tableaux, we will use the \textit{odd plactic ring} $\ZPl$, which is the unital ring defined by
\begin{equation}\begin{split}
\text{generators:}\qquad&A\\
\text{relations:}\qquad&yzx=-yxz\qquad\text{if }x<y\leq z\qquad(K'),\\
&xzy=-zxy\qquad\text{if }x\leq y<z\qquad(K'').
\end{split}\end{equation}
When we want to emphasize that the alphabet in question is $\lbrace1,2,\ldots,n\rbrace$, we will sometimes write $\ZPl_n$ instead of $\ZPl$.  The relations $(K'),(K'')$ are called \text{elementary Knuth transformations}.   We define a map from the set of semistandard Young tableaux with entries in the alphabet $A$ to the odd plactic ring by
\begin{equation}\begin{split}
\lbrace\text{SSYTs}\rbrace&\to\ZPl,\\
T&\mapsto w_r(T).
\end{split}\end{equation}
Since both the relations $(K'),(K'')$ are transpositions of letters with a minus sign, $\ZPl_n$ sits as an intermediate quotient between a free algebra and the skew polynomial ring,
\begin{equation*}
\Z\langle x_1,\ldots,x_n\rangle\twoheadrightarrow\ZPl_n\twoheadrightarrow\opol_n,
\end{equation*}
where the first map sends $\xt_i$ to $i$ and the second map sends $i$ to $\xt_i$.  If $w$ is a word in $\ZPl_n$, we will write $\xt^w$ for the image of $w$ in $\opol_n$.  In particular, a semistandard Young tableau $T$ is sent to $\xt^{w_r(T)}$.

The utility of the plactic ring is in large part due to the following remarkable theorem.
\begin{thm}[\cite{Fulton}, Section 2.1]\label{thm-word-tableau} Every word is equivalent, via relations $(K')$ and $(K'')$, to the row word $w_r(T)$ of a unique tableau $T$.\end{thm}
\noindent Thus the set of all Young tableaux with entries in $A$ forms a basis of $\ZPl$.  We will informally refer to the multiplication of tableaux in the following; what we mean is the multiplication of their row words in $\ZPl$.  In terms of tableaux, the relations $(K')$ and $(K'')$ can be interpreted as ``bumping transformations'':
\begin{equation*}\begin{split}
(K')\qquad&\young(yz)\cdot\young(x)=-\young(xz,y)\qquad\text{if }x<y\leq z,\\
(K'')\qquad&\young(xz)\cdot\young(y)=-\young(xy,z)\qquad\text{if }x\leq y<z.
\end{split}\end{equation*}
For a detailed exposition of bumping, see Section 1.1 of \cite{Fulton}.

\vspace{0.07in}

If a word $w$ is known to be the row word of some tableau, then it is easy to reconstruct the tableau from the word.  Since the row entries of a tableau never decrease from left to right and the column entries must always increase from top to bottom, reading the word $w$ from left to right until the first adjacent decreasing pair simply gives the bottom row of the tableau.  Then continuing to read until the next adjacent decreasing pair gives the second to bottom row, and so forth.  

\begin{example} Using $\Z_{>0}$ as the ordered alphabet,
\begin{equation*}
w=53422331112\qquad\text{corresponds to}\qquad\young(1112,2233,34,5).
\end{equation*}
\end{example}

\begin{defn} Let $\lambda$ be a partition.  Define an element of $\ZPl_n$ by
\begin{equation}
\widehat{s}_\lambda=(-1)^{dN(\lambda)+N(\lambda)}\sum_{T\in\SSYT(\lambda)}T.
\end{equation}
Its image in $\opol_n$,
\begin{equation}\label{eqn-defn-s-plactic}
s_\lambda^p=(-1)^{dN(\lambda)+N(\lambda)}\sum_{T\in\SSYT(\lambda)}\xt^{w_r(T)},
\end{equation}
is called the \textit{plactic odd Schur function} corresponding to $\lambda$.\end{defn}

The label $p$ in the notation for the plactic odd Schur function is to distinguish it from the combinatorial odd Schur functions $s_\lambda^K$ of \cite{EK} and the odd-symmetrized Schur functions $s_\lambda^s$ of \cite{EKL}; once we prove these three objects are equal later in this section, we will drop the extra labels.

\begin{rem} For the rest of the paper, we will not always specify whether we are working in $\osym$ or in $\osym_n$.  Generally speaking, results will hold in $\osym$.  The only change required in passing to $\osym_n$ is the understanding that certain elements become zero.  It is easy to see that $\osym_n=\osym/(e_m:m>n)$, and it will follow from Theorem \ref{thm-schur-coincidence} that $s_\lambda=0$ in $\osym_n$ if and only if $\lambda$ has height greater than $n$.  With this understood, the proofs and results in the rest of this paper work in either context.\end{rem}

\begin{lem}\label{lem-schur-straight} If $\lambda=(1^k)$, then up to sign, all three Schur functions coincide with the corresponding elementary polynomial:
\begin{equation}
s_{(1^k)}^K=s_{(1^k)}^p=s_{(1^k)}^s=(-1)^{\frac{1}{2}k(k-1)}e_k.
\end{equation}
If $\lambda=(k)$, the combinatorial and plactic Schur functions coincide with the corresponding complete polynomial:
\begin{equation}
s_{(k)}^K=s_{(k)}^p=h_k.
\end{equation}\end{lem}
The last equation equals $s_{(k)}^s$ too, but we will prove this later.
\begin{proof} The equality $s_{(1^k)}^K=(-1)^{\frac{1}{2}k(k-1)}e_k$ follows from Proposition 3.10 of \cite{EK} and the equality $s_{(1^k)}^s=(-1)^{\frac{1}{2}k(k-1)}e_k$ is Lemma 2.25 of \cite{EKL}.  Since the row word of a semistandard Young tableau on the shape $(1^k)$ is $i_k\cdots i_1$ for positive integers $i_1<\ldots<i_k$ and it takes $\frac{1}{2}k(k-1)$ transpositions to ascend-sort the monomial $\xt_{i_k}\cdots\xt_{i_1}$, the equality with $s_{(1^k)}^p$ holds.  It is similar but easier to show $s_{(k)}^p=h_k$, and $s_{(k)}^K=h_k$ is obvious because the only semistandard Young tableau of content $(k)$ is a row of 1's.\end{proof}

%
\subsection{Comparison with previous definitions}\label{subsec-comparison}
%

For a partition $\lambda$, let $\frac{i}{\lambda}$ be the Young diagram obtained by removing rows 1 through $i$ from the diagram corresponding to $\lambda$.  Similarly, let $\frac{\lambda}{i}$, $i|\lambda$, and $\lambda|i$ be obtained by removing rows $i$ through the bottom, columns $1$ through $i$, and columns $i$ through the rightmost respectively.  We say that a skew shape is a \textit{vertical strip} (respectively \textit{horizontal strip}) if no two of its boxes are in the same row (respectively column).  We say that a diagram $\mu$ is obtained from $\lambda$ by adding a vertical strip if $\lambda\subset\mu$ and $\mu/\lambda$ is a vertical strip; likewise for horizontal strips.  The following proposition was proved in \cite{EKL}.
\begin{prop}[Odd Pieri rule, $e$-right odd-symmetrized version]\label{prop-pieri-s} Let $\lambda$ be a partition.  Then
\begin{equation}\label{eqn-pieri-s}
s_\lambda^ss_{(1^k)}^s=\sum_\mu(-1)^{\left|\frac{i_1}{\lambda}\right|+\ldots+\left|\frac{i_k}{\lambda}\right|}s_\mu^s.
\end{equation}
The sum is over all $\mu$ obtained from $\lambda$ by adding a vertical strip of size $k$, and $i_1,\ldots,i_k$ are the rows of $\lambda$ to which a box was added.\end{prop}

The plactic odd Schur functions satisfy the same relation.  We prove the horizontal strip variant instead, for simplicity.
\begin{prop}[Odd Pieri rule, $h$-right plactic version]\label{prop-pieri-p} Let $\lambda$ be a partition.  Then
\begin{equation}\label{eqn-pieri-p}
(-1)^{NE(\lambda)}s_\lambda^ps_{(k)}^p=\sum_\mu(-1)^{NE(\mu)}(-1)^{\left|i_1|\lambda\right|+\ldots+\left|i_k|\lambda\right|}s_\mu^p.
\end{equation}
The sum is over all $\mu$ obtained from $\lambda$ by adding a horizontal strip of size $k$, and $i_1,\ldots,i_k$ are the columns of $\lambda$ to which a box was added.\end{prop}
\begin{proof} Expanding all Schur functions as odd plactic sums, we want to prove
\begin{equation*}
(-1)^{dN(\lambda)+N(\lambda)+NE(\lambda)}\xt^{w_r(T)}\xt^{w_r(V)}=(-1)^{dN(\mu)+N(\mu)+NE(\mu)}(-1)^{\left|i_1|\lambda\right|+\ldots+\left|i_k|\lambda\right|}\xt^{w_r(U)}
\end{equation*}
whenever $T$ is a Young tableau of shape $\lambda$, $V$ is a Young tableau of shape $(k)$, $TV=U$ in the even plactic ring with $\sh(U)=\mu$, and the boxes of $\mu/\lambda$ are in columns $i_1,\ldots,i_k$.  This is because mod 2 the even and odd plactic rings are isomorphic (by the obvious map, $w_r(T)\mapsto w_r(T)$), so the set of products $\xt^{w_r(T)}\xt^{w_r(V)}$ and the set of terms $\xt^{w_r(U)}$ which occur on the right-hand side of \eqref{eqn-pieri-p} are in bijection, with $(T,V)$ corresponding to $U$ if and only if $TV=U$ in the even plactic ring.  Suppose the leftmost box of $V$ has entry $j$.  If that box ends up in column $i$, we claim
\begin{equation*}
(-1)^{dN(\lambda)+N(\lambda)+NE(\lambda)}\xt^{w_r(T)}\xt_j=(-1)^{dN(\mu)+N(\mu)+NE(\mu)+\left|i|\lambda\right|}\xt^{w_r(U)}.
\end{equation*}
For any partition $\nu$, $(-1)^{dN(\nu)+N(\nu)+NE(\nu)}=(-1)^{NW(\nu)}$.  The new box of $\mu$ is not NorthWest of any other box (since it must be a southeast corner), so the sign discrepancy only counts those boxes NorthWest of the new box.  The sign $(-1)^{\left|i|\lambda\right|}$ counts boxes NorthEast of the new box, so the overall sign is $(-1)^{\sum_j(\lambda_j-1)}$.  And this is precisely the sign between $\xt^{w_r(T)}\xt_j$ and $\xt^{w_r(U)}$, since bumping a box past a row of length $r$ incurs a sign of $(-1)^{r-1}$.  Finally, note that as the boxes of $V$ are added one at a time, the signs cancel telescopically so as to yield the sign of equation \eqref{eqn-pieri-p}.\end{proof}

We now have the tools necessary to prove the main result of this section.
\begin{thm}\label{thm-schur-coincidence} The three notions of odd Schur function all coincide: for any partition $\lambda$,
\begin{equation}
s_\lambda^K=s_\lambda^p=s_\lambda^s.
\end{equation}\end{thm}
\begin{proof} The proof that $s_\lambda^K=s_\lambda^p$ is similar in spirit to the proof of the odd Pieri rule, $h$-right plactic version.  By the same sort of analysis as in that proof, one shows that for a partition $\mu=(\mu_1,\ldots,\mu_r)$,
\begin{equation*}\begin{split}
h_\mu=h_{\mu_1}\cdots h_{\mu_r}&=\sum_{\ct(T)=\mu}(-1)^{NE(T)+NE^<(T)+dN(T)+N(T)}\xt^{w_r(T)}\\
&=\sum_\lambda\sum_{T\in\SSYT(\lambda,\mu)}(-1)^{NE(\lambda)+NE^<(T)}s_\lambda^p.
\end{split}\end{equation*}
Since $h_\mu=\sum_\lambda K_{\lambda\mu}s_\lambda^K$ and both $\lbrace h_\mu\rbrace,\lbrace s_\lambda^K\rbrace$ are integral bases of the ring of odd symmetric functions, it suffices to check
\begin{equation*}
(-1)^{NE(\lambda)+NE^<(T)}=(-1)^{N(\lambda)+N^<(T)}
\end{equation*}
whenever $T$ is a Young tableau of shape $\lambda$ (the right-hand side is the sign with which $T$ is counted in defining $K_{\lambda\mu}$; see Section \ref{subsec-odd-symmetric}).  The sum computing the sign from a particular box $B$ of $T$ involves two types of other box: those Northwest and those NorthEast of $B$.  For those NorthEast, the sign is identical.  Those Northwest are ignored by the left-hand sign, but the entry of such a box is necessarily less than that of $B$ (since $T$ is semistandard), so the right-hand sign is $+1$.  Hence $s_\lambda^K=s_\lambda^p$.

We now use the Pieri rules (Propositions \ref{prop-pieri-s}, \ref{prop-pieri-p}) to prove $s_\lambda^p=s_\lambda^s$.  For $\lambda=(1^k)$, this is true by Lemma \ref{lem-schur-straight}.  Using this as a base case, we now induct on the width of $\lambda$, and within each particular width we induct with respect to the lexicographic order of $\lambda^T$.  Since we have shown $s_\lambda^K=s_\lambda^p$, we can apply the involution $\psi_1\psi_2$ to equation \eqref{eqn-pieri-p} to show the plactic Schur functions obey an $e$-right Pieri rule of the same form and signs as the one obeyed by the Schur functions $s_\lambda^K$,
\begin{equation}\label{eqn-pieri-p-2}
s_\lambda^ps_{(1^k)}^p=\sum_\mu(-1)^{\left|\frac{i_1}{\lambda}\right|+\ldots+\left|\frac{i_k}{\lambda}\right|}s_\mu^p.
\end{equation}
Let $\lambda$ be a partition of width $\lambda_1=r$.  Using the $e$-right Pieri rule, both types of Schur function satisfy
\begin{equation*}
s_{(\lambda^T_1)}\cdots s_{(\lambda^T_r)}=\sum_\mu\pm s_\mu,
\end{equation*}
where each $\mu$ has width at most $r$ and is lexicographically greater than or equal to $\lambda$.  The coefficient of $s_\lambda$ on the right-hand side is $\pm1$.  All the signs $\pm$ are the same for the two types of Schur function, so this allows us to solve for both $s_\lambda^p,s_\lambda^s$ in terms of elementary functions by the same expressions; hence $s_\lambda^p=s_\lambda^s$.\end{proof}

\begin{cor} The span of the $\widehat{s}_\lambda$ in $\ZPl_n$ is a subalgebra isomorphic to $\osym_n$, and this subalgebra is taken isomorphically onto $\osym_n\subset\opol_n$ by the map $w\mapsto\xt^w$.\end{cor}

For the rest of this paper, we will drop the superscript labels on odd Schur functions.

%
\section{Littlewood-Richardson rules}
%

%
\subsection{The even Littlewood-Richardson rule}\label{subsec-even-lr}
%

\quad\\
\begin{center}\textit{For this section only, we work in the even setting.}\end{center}
\quad\\
Let $\mu,\nu,\lambda$ be partitions.  The \textit{Littlewood-Richardson coefficient} $c_{\mu\nu}^\lambda$ is the coefficient of $s_\lambda$ in $s_\mu s_\nu$,
\begin{equation*}
s_\mu s_\nu=\sum_\lambda c_{\mu\nu}^\lambda s_\lambda.
\end{equation*}
If $|\mu|+|\nu|\neq|\lambda|$, then $c_{\mu\nu}^\lambda=0$.  Schur functions are generating functions for semistandard Young tableaux of a given shape,
\begin{equation*}
s_\lambda=\sum_{T\in\SSYT(\lambda)}x^{w_r(T)}.
\end{equation*}
It follows that, for any fixed semistandard Young tableau $T_0$ of shape $\lambda$,
\begin{equation}\label{eqn-easy-even-lr}
c_{\mu\nu}^\lambda=\#\lbrace U\in\SSYT(\mu),V\in\SSYT(\nu):UV=T_0\rbrace.
\end{equation}
Here, the product of tableaux is taken in the (even) plactic ring.  If $T_0=T_\lambda$, then it is not hard to see that $UV=T_\lambda$ forces $V$ to be $T_\nu$; equation \eqref{eqn-easy-even-lr} becomes
\begin{equation}\label{eqn-easy-even-lr-2}
c_{\mu\nu}^\lambda=\#\lbrace U\in\SSYT(\mu):UT_\nu=T_\lambda\rbrace.
\end{equation}
Since computing the set on the right-hand side of \eqref{eqn-easy-even-lr-2} can be tricky in practice, we would like to have a simpler combinatorial description of the $c_{\mu\nu}^\lambda$.  The Littlewood-Richardson rule provides one of many such simpler descriptions.  Good accounts of the rule and its proof are given in \cite[Chapter 5]{Fulton} and \cite[Appendix A1]{Stanley2}.  We will review the terminology and statement here.

\vspace{0.07in}

Recall that a \textit{(Young) skew shape} $\lambda/\mu$ is the complement of a subdiagram $\mu\subseteq\lambda$.  A \textit{semistandard skew tableau} is a skew shape which has been filled with entries from some ordered alphabet, subject to the same rules as for a semistandard tableau: entries must strictly increase in columns (top to bottom) and must not decrease in rows (left to right).  We write $\SSYT(\lambda/\mu)$ for the set of semistandard skew tableaux of shape $\lambda/\mu$ and $\SSYT(\lambda/\mu,\nu)$ for the set of semistandard skew tableaux of shape $\lambda/\mu$ and content $\nu$.

A word $w=w_1\cdots w_r$ in some ordered alphabet is called \textit{Yamanouchi} (or a \textit{reverse lattice word}) if, when read backwards, each initial subword has at least as many $a$'s as $b$'s whenever $a<b$.  For example, $312211$ is Yamanouchi but $1221$ and $112$ are not.  A skew tableau $S$ is called a \textit{Littlewood-Richardson tableau} if $w_r(S)$ is a Yamanouchi word.  The following is Proposition 3, Chapter 5 of \cite{Fulton} and Theorem A1.3.3 of \cite{Stanley2}.
\begin{thm}[(Even) Littlewood-Richardson Rule]\label{thm-even-lr} The coefficient $c_{\mu\nu}^\lambda$ equals the number of Littlewood-Richardson tableaux $S$ of shape $\lambda/\mu$ and content $\nu$.\end{thm}
One specific bijection between the set described in the theorem and the set in equation \eqref{eqn-easy-even-lr-2} is described in the following section, where we use it to deduce an odd analogue of Theorem \ref{thm-even-lr}.

\begin{example} Let $\mu\subseteq\lambda$ be a subdiagram and let $k=|\lambda|-|\mu|\geq1$.  If $S$ is a Littlewood-Richardson tableau of shape $\lambda/\mu$ and content $(k)$, then no two boxes of $S$ can be in the same column (since all entries equal 1).  And on any such skew shape $\lambda/\mu$, there is exactly one tableau of content $(k)$.  Thus $c_{\mu(k)}^\lambda=1$ if $\lambda/\mu$ is a horizontal strip and equals 0 otherwise.  We have deduced the (even) Pieri rule,
\begin{equation}
s_\mu s_{(k)}=\sum_{\substack{\lambda\\\lambda/\mu\text{ is a}\\\text{horizontal strip}}}s_\lambda.
\end{equation}
Using the standard involution $\omega$ on $\sym$, or just arguing in analogy with the above, the same is true if $(k)$ is replaced by $(1^k)$ and ``horizontal'' is replaced by ``vertical.''\end{example}

\begin{example} The lowest degree product which is not described by the Pieri rule is 
\begin{equation*}
s_{21}s_{21}=s_{2211}+s_{222}+s_{3111}+2s_{321}+s_{33}+s_{411}+s_{42}.
\end{equation*}
\end{example}

%
\subsection{The odd Littlewood-Richardson rule}\label{subsec-odd-lr}
%

The sign of a Young tableau $T$ is the sign between its row word monomial $\xt^{w_r(T)}$ and the ascend-sorting of that monomial.  As explained in Section \ref{subsec-odd-symmetric}, this sign equals $(-1)^{N^<(T)}$.  If $S$ is a skew tableau of shape $\lambda/\mu$, let $j$ be an element of the alphabet less than every entry of $S$ and let $\widehat{S}$ be the Young tableau of shape $\lambda$ formed by placing $j$ in each box of $\mu$ and filling the rest so as to match $S$.  We then define the sign of $S$ to be
\begin{equation*}
\sign(S)=\sign(\widehat{S})=N^<(\widehat{S}).
\end{equation*}
When $\mu=(0)$, this reduces to the sign of a tableau as defined earlier.  More generally, whenever we write a count $E, nW^\geq, dS,\ldots$ evaluated on $S$, read $\widehat{S}$ for $S$.

\begin{example} To either alphabet $\Z_{>0}$ or $\lbrace1,2,\ldots,n\rbrace$, we can always adjoin 0 and take $j=0$.  If $\lambda=(3,3,2,1)$, $\mu=(2,1,1)$, and
\begin{equation*}
S=\young(::1,:12,:2,3),\quad\text{then}\quad\widehat{S}=\young(001,012,02,3)\quad\text{and}\quad\sign(S)=(-1)^{18}=1.
\end{equation*}\end{example}

\begin{rem} We consider the partition $\mu$ to be part of the data of $S$.  For example, both $(1,1)/(1)$ and $(2)/(1)$ consist of a single box, but if we fill these two boxes with equal entries, the resulting skew tableaux have opposite signs.\end{rem}

\begin{defn} Let $\mu,\nu,\lambda$ be partitions.  The \textit{odd Littlewood-Richardson coefficient} $c_{\mu\nu}^\lambda$ is the coefficient of $s_\lambda$ when $s_\mu s_\nu$ is expanded in the basis of odd Schur functions,
\begin{equation*}
s_\mu s_\nu=\sum_\lambda c_{\mu\nu}^\lambda s_\lambda.
\end{equation*}
If $|\mu|+|\nu|\neq|\lambda|$, then $c_{\mu\nu}^\lambda=0$.\end{defn}

Note that the odd Pieri rules compute certain odd Littlewood-Richardson coefficients.  Using the involution $\psi_1\psi_2$ and the anti-involution $\psi_3$, we know the odd Littlewood-Richardson coefficient $c_{\mu\nu}^\lambda$ whenever $\mu$ or $\nu$ has either height or width 1.

If $Y,Z$ are two nonzero monomials in $\opol_n$ such that $Y=\pm Z$, then let $\sign(Y,Z)$ denote the sign between them; for example $\sign(x_1x_2x_3,x_1x_3x_2)=-1$.
\begin{lem}\label{lem-easy-lr} The odd Littlewood-Richardson coefficient $c_{\mu\nu}^\lambda$ is
\begin{equation}\label{eqn-easy-lr}
c_{\mu\nu}^\lambda=(-1)^{dN(\mu)+dN(\nu)+dN(\lambda)+N(\mu)+\sum_i(\lambda/\nu)_i\left|\frac{\nu}{i}\right|}\sum_{\substack{U\in\SSYT(\mu)\\UT_\nu=T_\lambda}}(-1)^{N^<(U)}.
\end{equation}
In the summation condition, the product $UT_\nu$ is taken in the even plactic ring (so it is an equality up to sign in $\ZPl_n$).\end{lem}
\begin{proof} First note that
\begin{equation}\begin{split}
c_{\mu\nu}^\lambda&=(-1)^{dN(\mu)+dN(\nu)+dN(\lambda)+N(\mu)+N(\nu)+N(\lambda)}\sum_{\substack{U\in\SSYT(\mu)\\V\in\SSYT(\nu)\\UV=T_\lambda}}\sign(\xt^{w_r(U)}\xt^{w_r(V)},\xt^{w_r(T_\lambda)})\\
&=(-1)^{dN(\mu)+dN(\nu)+dN(\lambda)+N(\mu)+N(\nu)+N(\lambda)}\sum_{\substack{U\in\SSYT(\mu)\\UT_\nu=T_\lambda}}\sign(\xt^{w_r(U)}\xt^{w_r(T_\nu)},\xt^{w_r(T_\lambda)}).
\end{split}\end{equation}
The first equality is immediate from the definition of $c_{\mu\nu}^\lambda$ and equation \eqref{eqn-defn-s-plactic}.  The second follows from the following fact: if $U,V$ are semistandard tableaux of shapes $\mu,\nu$ and $UV=T_\lambda$ in the even plactic ring, then $V=T_\nu$ \cite[Section 5.2]{Fulton}.  To turn $\xt^{w_r(U)}\xt^{w_r(T_\nu)}$ into $\xt^{w_r(T_\lambda)}$, we proceed in two steps:
\begin{enumerate}
\item Ascend-sort each of $\xt^{w_r(U)}$, $\xt^{w_r(T_\nu)}$, and $\xt^{w_r(T_\lambda)}$ separately.  This incurs the sign $(-1)^{N^<(U)+N(\nu)+N(\lambda)}$.  The monomials are now $:\xt^{w_r(T_{\lambda/\nu})}::\xt^{w_r(T_\nu)}:$ and $:\xt^{w_r(T_\lambda)}:$, where the colons denote ascend-sorting.
\item To sort these monomials together requires $\sum_i(\lambda/\nu)_i\left|\frac{\nu}{i}\right|$ transpositions.
\end{enumerate}
The lemma follows.\end{proof}

\begin{thm}[Odd Littlewood-Richardson rule]\label{thm-lr} The odd Littlewood-Richardson coefficient $c_{\mu\nu}^\lambda$ is
\begin{equation}\label{eqn-lr}
c_{\mu\nu}^\lambda=(-1)^{N(\mu)+N(\lambda)}\sum_{\substack{S\in\SSYT(\lambda/\mu,\nu)\\w_r(S)\text{ is Yamanouchi}}}(-1)^{N^<(S)}.
\end{equation}\end{thm}
\begin{proof} The Littlewood-Richardson bijection of \cite[Section 5.2]{Fulton},
\begin{equation}\label{eqn-lr-bij}
\lbrace S\in\SSYT(\lambda/\mu,\nu):w_r(S)\text{ is Yamanouchi}\rbrace\to\lbrace U\in\SSYT(\mu,\lambda/\nu):\\UT_\nu=T_\lambda\rbrace,
\end{equation}
is described as follows.  Let $U_0$ be a semistandard Young tableau of shape $\mu$, all of whose entries are less than all of those of $S$.  Then form $T_S\in\SSYT(\lambda)$ by giving it the entries of $U_0$ on $\mu\subseteq\lambda$ and the entries of $S$ on $\lambda/\mu$.  Under the RSK correspondence, 
\begin{equation*}
(T_\lambda,T_S)\leftrightarrow\left(\substack{\underline{t}\hspace{0.02in}\underline{u}\\\underline{x}\hspace{0.02in}\underline{v}}\right),
\end{equation*}
where $\underline{x},\underline{t}$ have length $|\mu|$ and $\underline{u},\underline{v}$ have length $|\nu|$.  We can describe what these sub-words correspond to under the RSK correspondence as well:
\begin{equation*}
(U,U_0)\leftrightarrow\left(\substack{\underline{t}\\\underline{x}}\right),\qquad(T_\nu,T_\nu)\leftrightarrow\left(\substack{\underline{u}\\\underline{v}}\right)
\end{equation*}
for some $U\in\SSYT(\mu,\lambda/\nu)$.  The Littlewood-Richardson bijection \eqref{eqn-lr-bij} assigns $U$ to $S$.

In order to prove the theorem, we have to relate the quantities $N^<(U)$ and $N^<(S)$.  The odd RSK correspondence of \cite{EK} allows us to do this; it implies
\begin{equation}\label{eqn-pf-lr-1}\begin{split}
\sign(\underline{x}\hspace{0.02in}\underline{v})&=(-1)^{dN(\lambda)+N(\lambda)+N^<(T_S)}\\
&=(-1)^{dN(\lambda)+N(\lambda)+N^<(U_0)+N^<(S)},\\
\sign(\underline{x})&=(-1)^{dN(\mu)+N^<(U)+N^<(U_0)},\\
\sign(\underline{v})&=(-1)^{dN(\nu)}.
\end{split}\end{equation}
Since we know the contents of $\underline{x}$, $\underline{v}$, and $\underline{x}\hspace{0.02in}\underline{v}$, the sign of the latest can be expressed in terms of the former two,
\begin{equation}\label{eqn-pf-lr-2}
\sign(\underline{x}\hspace{0.02in}\underline{v})=\sign(\underline{x})\sign(\underline{v})(-1)^{\sum_i(\lambda/\nu)_i\left|\frac{\nu}{i}\right|}.
\end{equation}
Comparing the signs in equations \eqref{eqn-pf-lr-1} and \eqref{eqn-pf-lr-2}, 
\begin{equation*}
(-1)^{dN(\mu)+dN(\nu)+dN(\lambda)+N(\lambda)+N^<(U)+N^<(S)+\sum_i(\lambda/\nu)_i\left|\frac{\nu}{i}\right|}=1.
\end{equation*}
Applying this to Lemma \ref{lem-easy-lr}, the theorem follows.\end{proof}

\begin{example} The first interesting cancellation is $c_{(2,1)(2,1)}^{(3,2,1)}=0$ (in the even case, it equals 2).  The Littlewood-Richardson skew tableaux of shape $(3,2,1)/(2,1)$ and content $(2,1)$ are
\begin{equation*}
\young(::1,:1,2),\quad\young(::1,:2,1).
\end{equation*}
They have $N^<(S)$ equal to $7$ and $6$, respectively.\end{example}

\begin{rem} It follows from equation \eqref{eqn-psi-schur} that
\begin{equation}\begin{split}
&c_{\mu\nu}^\lambda=(-1)^{dN(\mu)+dN(\nu)+dN(\lambda)+N(\mu)+N(\nu)+N(\lambda)}c_{\nu\mu}^\lambda,\\
&c_{\mu\nu}^\lambda=(-1)^{NE(\mu)+NE(\nu)+NE(\lambda)}c_{\mu^T\nu^T}^{\lambda^T}.
\end{split}\end{equation}
These symmetries constrain some signs associated to Young diagrams.  For instance, if $c_{\mu\mu}^\lambda\neq0$, then $(-1)^{dN(\lambda)+N(\lambda)}=1$.
\end{rem}

%
\subsection{Translation to Knutson-Tao hives}\label{subsec-hives}
%

There are many combinatorial expressions for the even Littlewood-Richardson coefficients.  Among them, the several which are expressible in terms of integer points of rational convex polytopes are especially interesting; one reason is that as result of Knutson and Tao's proof of the Saturation Conjecture \cite{KnutsonTao}, Klyachko's system of inequalities \cite{Klyachko} gives a necessary and sufficient criterion for a Littlewood-Richardson coefficient to be nonzero.  

These expressions in terms of polytopes include Gelfand-Zeitlin (GZ) patterns \cite{GelfandZelevinsky}, Berenstein-Zelevinsky triangles \cite{BerensteinZelevinsky}, the Littlewood-Richardson triangles of Pak and Vallejo \cite{PakVallejo}, and the honeycombs and hives of Knutson and Tao \cite{KnutsonTao}.  As explained in the exposition of \cite{Buch}, the hive model is particularly convenient and flexible.  In this section, we will write down the bijection of \cite{PakVallejo} from Littlewood-Richardson skew tableaux to hives, by way of Littlewood-Richardson triangles.  In both the triangle and hive settings, the sign associated to a diagram will come from a quadratic form on the ambient space of the polytope in question.

\vspace{0.07in}

The following exposition follows \cite{PakVallejo} rather closely; our only contribution is to track the signs which arise in the odd case.  For fixed $n\geq1$, we will work in $V_\Z=\Z^{\binom{n+2}{2}}$ and $V=V_\Z\otimes_\Z\R$, with coordinates $\lbrace a_{i,j}:0\leq i\leq j\leq n\rbrace$.  We call elements of $V$ \textit{triangles} and write them pictorially as (for $n=3$)
\begin{equation*}
\begin{matrix}&&&a_{0,0}&&&\\\\&&a_{0,1}&&a_{1,1}&&&\\\\&a_{0,2}&&a_{1,2}&&a_{2,2}&&\\\\a_{0,3}&&a_{1,3}&&a_{2,3}&&a_{3,3}.\end{matrix}
\end{equation*}

\begin{defn} A \textit{Littlewood-Richardson triangle} is an element $A=(a_{i,j})\in V$ such that
\begin{enumerate}
\item $a_{i,j}\geq0$ for $1\leq i<j\leq n$,
\item
\begin{equation*}
\sum_{p=0}^{i-1}a_{p,j}\geq\sum_{p=0}^ia_{p,j+1}\qquad\text{for }1\leq i\leq j<n,\text{ and}
\end{equation*}
\item
\begin{equation*}
\sum_{q=i}^ja_{i,q}\geq\sum_{q=i+1}^{j+1}a_{i+1,q}\qquad\text{for }1\leq i<j<n.
\end{equation*}
\end{enumerate}
Note that $a_{0,j}<0$ and $a_{j,j}<0$ are both possible.\end{defn}
Let $\bigtriangleup_{LR}$ denote the set of all Littlewood-Richardson triangles, a cone in $V$.  Given a Littlewood-Richardson triangle $A=(a_{i,j})\in\bigtriangleup_{LR}$, define partitions $\lambda,\mu,\nu$ by
\begin{equation*}
\lambda_j=\sum_{p=0}^ja_{p,j},\qquad\mu_j=a_{0,j},\qquad\nu_i=\sum_{q=i}^ka_{i,q}.
\end{equation*}
Let $\bigtriangleup_{LR}(\lambda,\mu,\nu)$ be the set of Littlewood-Richardson triangles with fixed $\lambda,\mu,\nu$.  Each set $\bigtriangleup_{LR}(\lambda,\mu,\nu)$ is a convex polytope in $V$.

To a Littlewood-Richardson skew tableau $S\in\SSYT(\lambda/\mu,\nu)$, associate an element $A_S$ of $\bigtriangleup_{LR}(\lambda,\mu,\nu)$ by setting
\begin{equation}\begin{split}
&a_{0,0}=0,\qquad a_{0,j}=\mu_j,\\
&a_{i,j}=\#\lbrace\text{entries equal to }i\text{ in row }j\text{ of }S\rbrace\quad(0<i\leq j).
\end{split}\end{equation}
\begin{lem}[\cite{PakVallejo}, Lemma 3.1] Suppose $|\mu|+|\nu|=|\lambda|$.  Then the assignment $S\mapsto A_S$ is a bijection between the set of Littlewood-Richardson tableaux in $\SSYT(\lambda/\mu,\nu)$ and the set of integer points of $\bigtriangleup_{LR}(\lambda,\mu,\nu)$.\end{lem}
It is easy to read off $N^<(S)$ from the triangle $A_S=(a_{i,j})$.  For each entry $a_{i,j}$, write $Y_{i,j}$ for the sum of all the $a_{p,q}$ in the shaded region below, where the $(i,j)$ place is the dot drawn with a hollow center:
\vspace{0.15in}
\begin{center}\includegraphics[width=2.5in]{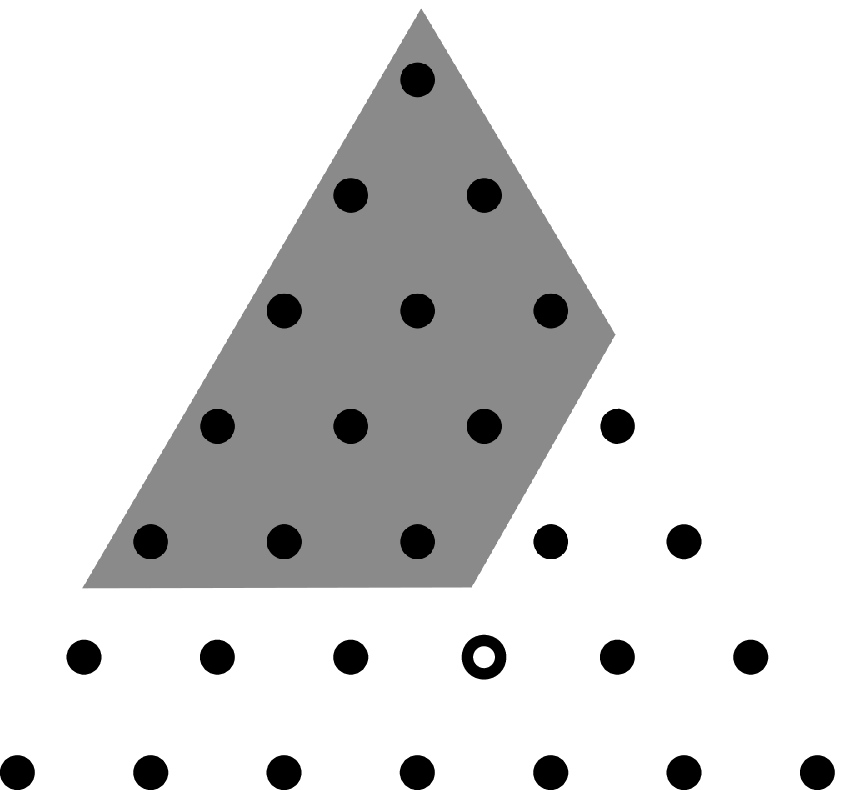}\end{center}
\vspace{0.15in}
More formally,
\begin{equation}
Y_{i,j}=\sum_{p=0}^{i-1}\sum_{q=p}^{j-1}a_{p,q}.
\end{equation}
Consider the quadratic form
\begin{equation}
Q_\bigtriangleup(A)=\sum_{i,j}a_{i,j}Y_{i,j}=\sum_{i=0}^n\sum_{j=i}^n\left(a_{i,j}\sum_{p=0}^{i-1}\sum_{q=p}^{j-1}a_{p,q}\right).
\end{equation}
Then it is immediate from the description of the bijection above that
\begin{equation*}
N^<(S)=Q_\bigtriangleup(A_S),
\end{equation*}
so
\begin{equation}
c_{\mu\nu}^\lambda=(-1)^{N(\mu)+N(\lambda)}\sum_{A\in\bigtriangleup_{LR}(\lambda,\mu,\nu)\cap V_\Z}(-1)^{Q_\bigtriangleup(A)},
\end{equation}
as long as $\lambda,\mu,\nu$ all have at most $n$ parts.

\begin{example}\label{ex-triangle} If
\begin{equation*}
S=\young(::1,:2,1),
\end{equation*}
then
\begin{equation*}
A_S=\begin{matrix}&&&0&&&\\\\&&2&&1&&\\\\&1&&0&&1&&\\\\0&&1&&0&&0&\end{matrix}
\end{equation*}
and $Q_\bigtriangleup(A_S)=6$.
\end{example}

We now translate this result into the language of hives.
\begin{defn} A \textit{hive} is an element $H=(h_{i,j})\in V$ with $h_{0,0}=0$ which satisfies the inequalities
\begin{equation}\label{eqn-hive-inequalities}\begin{split}
&(R)\qquad h_{i,j}-h_{i,j-1}\geq h_{i-1,j}-h_{i-1,j-1}\text{ for }1\leq i<j\leq n,\\
&(V)\qquad h_{i-1,j}-h_{i-1,j-1}\geq h_{i,j+1}-h_{i,j}\text{ for }1\leq i\leq j<n,\\
&(L)\qquad h_{i,j}-h_{i-1,j}\geq h_{i+1,j+1}-h_{i,j+1}\text{ for }1\leq i\leq j<n.
\end{split}\end{equation}\end{defn}
Let $\mathfrak{H}$ be the set of all hives; this is a cone in $V$.  The inequalities \eqref{eqn-hive-inequalities} have a geometric interpretation when we express $H$ as a triangle.  Inside a triangle, there are three types of rhombi which can be made out of two adjacent smallest-size triangles:
\vspace{0.15in}
\begin{center}\includegraphics[width=2.5in]{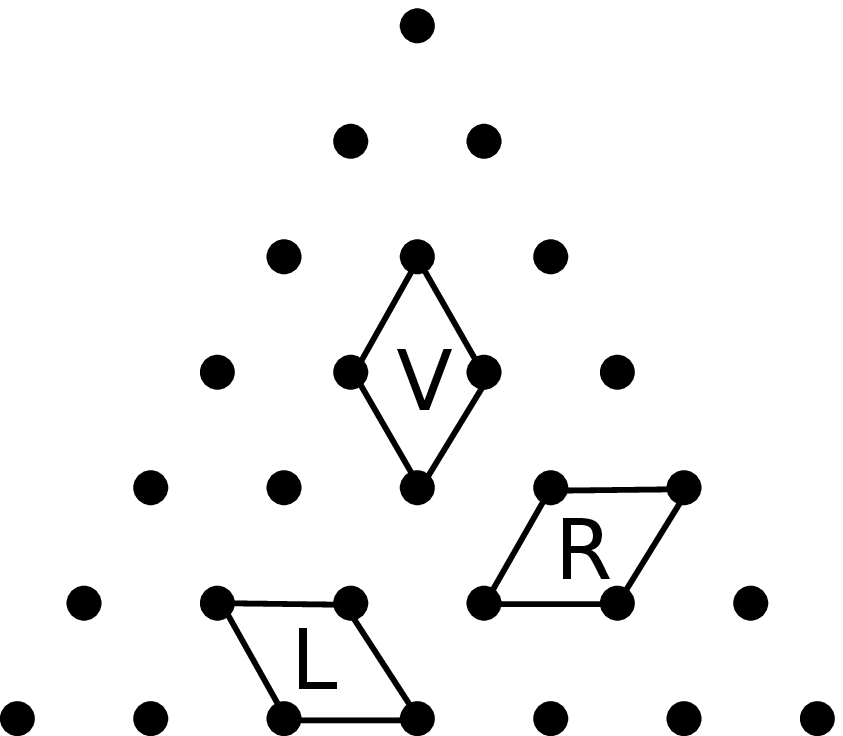}\end{center}
\vspace{0.15in}
We call these right-slanted (R), vertical (V), and left-slanted (L) rhombi.  The inequalities \eqref{eqn-hive-inequalities} say that the sum of the entries at the obtuse angles of any such rhombus is greater than or equal to the sum at the acute angles; the three inequalities are the right-slanted, vertical, and left-slanted cases, respectively.

As with Littlewood-Richardson triangles, we associate three partitions to a hive,
\begin{equation*}
\lambda_j=h_{j,j}-h_{j-1,j-1},\qquad\mu_j=h_{0,j}-h_{0,j-1},\qquad\nu_i=h_{i,n}-h_{i-1,n}.
\end{equation*}
Let $\mathfrak{H}(\lambda,\mu,\nu)$ be the set of all hives with corresponding partitions $\lambda,\mu,\nu$.  Each $\mathfrak{H}(\lambda,\mu,\nu)$ is a convex polytope in $V$.

\begin{thm}[\cite{PakVallejo}, Theorem 4.1] Let $\Phi:V\to V$ be the linear map which takes $A=(a_{i,j})$ to $H=(h_{i,j})$, where
\begin{equation}\label{eqn-a-to-h}
h_{i,j}=\sum_{p=0}^i\sum_{q=p}^ja_{p,q}.
\end{equation}
Then $\Phi$ is a volume-preserving isomorphism and induces bijections $\bigtriangleup_{LR}(\lambda,\mu,\nu)\to\mathfrak{H}(\lambda,\mu,\nu)$ for all $\lambda,\mu,\nu$.\end{thm}

As a matter of convention, let $h_{i,j}=0$ if either $i>j$, $i<0$, or $j>n$.  It follows from equation \eqref{eqn-a-to-h} and an inclusion-exclusion argument that
\begin{equation}\label{eqn-hij}\begin{split}
&h_{i,j}=a_{i,j}+h_{i-1,j}+h_{i,j-1}-h_{i-1,j-1}\text{ if }0\leq i<j\leq n,\\
&h_{i,i}=a_{i,i}+h_{i-1,i}\text{ if }0\leq i=j\leq n.
\end{split}\end{equation}
Let $Q_\mathfrak{H}=Q_\bigtriangleup\circ\Phi^{-1}$.  Then for a hive $H=(h_{i,j})$, equations \eqref{eqn-hij} imply
\begin{equation}
Q_\mathfrak{H}(H)=\sum_{i=1}^n\sum_{j=i}^nh_{i-1,j-1}\left(h_{i,j}-h_{i-1,j}-h_{i,j-1}+h_{i-1,j-1}\right)-\sum_{i=1}^{n-1}h_{i,i}^2.
\end{equation}
Note that the parenthesized term is non-negative by the right-slanted rhombus inequality $(R)$.  It follows that if $\mu,\nu,\lambda$ are partitions with at most $n$ parts, then
\begin{equation}
N^<(S)=Q_\bigtriangleup(A_S)=Q_\mathfrak{H}(\Phi(A_S)),
\end{equation}
so
\begin{equation}
c_{\mu\nu}^\lambda=(-1)^{N(\mu)+N(\lambda)}\sum_{H\in\mathfrak{H}(\lambda,\mu,\nu)\cap V_\Z}(-1)^{Q_\mathfrak{H}(H)}.
\end{equation}
The quantity $N(\mu)+N(\lambda)$ can also be expressed as a quadratic form in either the variables $(a_{i,j})$ or the variables $(h_{i,j})$.

\begin{example} With $S$ and $A_S$ as in Example \ref{ex-triangle},
\begin{equation}
\Phi(A_S)=\begin{matrix}&&&0&&&\\\\&&2&&3&&&\\\\&3&&4&&5&&\\\\3&&5&&6&&6&\end{matrix}
\end{equation}
and $Q_\mathfrak{H}(\Phi(A_S))=6$.
\end{example}


\bibliographystyle{alpha}
\bibliography{ellis-bib}

\begin{thebibliography}{Web10b}

\bibitem[Buc00]{Buch}
A.~Buch.
\newblock The saturation conjecture (after {A}. {K}nutson and {T}. {T}ao).
\newblock {\em L'Enseignement Math\'{e}matique, IIe S\'{e}rie}, 46:43--60,
  2000.
\newblock \href{http://arxiv.org/abs/math/9810180}{arXiv:math.CO/9810180}.

\bibitem[BZ92]{BerensteinZelevinsky}
A.~Berenstein and A.~Zelevinsky.
\newblock Triple multiplicities for $sl(r+1)$ and the spectrum of the exterior
  algebra of the adjoint representation.
\newblock {\em Journal of Algebraic Combinatorics}, 1:7--22, 1992.

\bibitem[EK11]{EK}
A.~P. Ellis and M.~Khovanov.
\newblock The {H}opf algebra of odd symmetric functions.
\newblock 2011.
\newblock \href{http://arxiv.org/abs/1107.5610}{arXiv:math.QA/1107.5610}.

\bibitem[EKL11]{EKL}
A.~P. Ellis, M.~Khovanov, and A.~Lauda.
\newblock The odd nil{H}ecke algebra and its diagrammatics.
\newblock 2011.
\newblock \href{http://arxiv.org/abs/1111.1320}{arXiv:math.QA/1111.1320}.

\bibitem[Ful97]{Fulton}
W.~Fulton.
\newblock {\em Young tableaux: With applications to representation theory and
  geometry}, volume~35 of {\em London Mathematical Society Student Texts}.
\newblock Cambridge University Press, Cambridge, 1997.

\bibitem[GZ86]{GelfandZelevinsky}
I.~Gelfand and A.~Zelevinsky.
\newblock Multiplicities and good bases for $gl_n$.
\newblock In {\em Group theoretical methods in physics, vol. {II} (Yurmala,
  1985)}, pages 147--159, Utrecht, 1986. VNU Scientific Press.
\newblock
  \href{http://www.math.neu.edu/zelevinsky/Gelfand-canonical-basis-gln.pdf}{ht%
tp://www.math.neu.edu/zelevinsky/Gelfand-canonical-basis-gln.pdf}.

\bibitem[Kly98]{Klyachko}
A.~Klyachko.
\newblock Stable bundles, representation theory and {H}ermitian operators.
\newblock {\em Selecta Mathematica, New Series}, 4(3):419--445, 1998.

\bibitem[KT99]{KnutsonTao}
A.~Knutson and T.~Tao.
\newblock The honeycomb model of {$GL_n(\C)$} tensor products {I}: Proof of the
  saturation conjecture.
\newblock {\em Journal of the American Mathematical Society}, 12(4):1055--1090,
  1999.
\newblock
  \href{http://www.ams.org/journals/jams/1999-12-04/S0894-0347-99-00299-4/}{ht%
tp://www.ams.org/journals/jams/1999-12-04/S0894-0347-99-00299-4/}.

\bibitem[ORS07]{ORS}
P.~Ozsv\'{a}th, J.~Rasmussen, and Z.~Szab\'{o}.
\newblock Odd {K}hovanov homology.
\newblock 2007.
\newblock \href{http://arxiv.org/abs/0710.4300}{arXiv:math.QA/0710.4300}.

\bibitem[PV05]{PakVallejo}
I.~Pak and E.~Vallejo.
\newblock Combinatorics and geometry of {L}ittlewood-{R}ichardson cones.
\newblock {\em European Journal of Combinatorics}, 26:995--1008, 2005.
\newblock
  \href{http://www.math.ucla.edu/~pak/papers/liri91.pdf}{http://www.math.ucla.%
edu/$\sim$pak/papers/liri91.pdf}.

\bibitem[Sta99]{Stanley2}
R.~P. Stanley.
\newblock {\em Enumerative Combinatorics, vol. 2}.
\newblock Cambridge University Press, Cambridge, UK, 1999.

\bibitem[Web10a]{Web}
B.~Webster.
\newblock Knot invariants and higher representation theory {I}: {D}iagrammatic
  and geometric categorification of tensor products, 2010.
\newblock \href{http://arxiv.org/abs/1001.2020}{arXiv:math.QA/1001.2020}.

\bibitem[Web10b]{Web2}
B.~Webster.
\newblock Knot invariants and higher representation theory {II}: {T}he
  categorification of quantum knot invariants, 2010.
\newblock \href{http://arxiv.org/abs/1005.4559}{arXiv:math.QA/1005.4559}.

\end{thebibliography}

%

\vspace{0.1in}

\noindent A.P.E.: { \sl \small Department of Mathematics, Columbia University, New
York, NY 10027, USA} \newline \noindent {\tt \small email: ellis@math.columbia.edu}

%
\end{document}